\documentstyle[a4,11pt,amssymb]{article}
\newenvironment{proof}{{\sc Proof of}}{\raisebox{.8ex}{\framebox[2mm][b]{}}}
\newtheorem{theorem}{\large\bf Theorem}
\newtheorem{proposition}{\large\bf Proposition}
\newtheorem{lemma}{\large\bf Lemma}

\newtheorem{corollary}{\large\bf Corollary}

\def\today{\number\day\space\ifcase\month\or
January\or February\or March\or April\or May\or June\or
July\or August\or September\or October\or November\or December\fi
\space\number\year}
\begin{document}
\title{ On limit theorems for continued fractions }

\author {Zbigniew S. Szewczak
\footnote{Nicolaus Copernicus
University, Faculty of Mathematics and Computer Science,
ul. Chopina 12/18, 87-100 Toru\'n, Poland, e-mail:
zssz@mat.uni.torun.pl}}
\maketitle

\begin{abstract}
It is shown that for sums of functionals of digits in continued fraction
expansion the Kolmogorov-Feller
weak laws of large numbers and
the Khinchine-L\'evy-Feller-Raikov characterization
of the domain of attraction of the normal law hold.

\end{abstract}

\section{ Introduction and result}

\noindent

Let $a_n(x),$ $n\in{\mathbb N}=\{1,2,\ldots\},$
denote the partial quotients (or digits)
in the simple non-terminating continued fraction expansion of an irrational
number $x\in (0,1],$ (Cf. \cite{Kn64})
\[
x =
{{1}\over{a_1(x) + {{1}\over{a_2(x) + {{1}\over{a_3(x)+\ldots}}}}}}=
{{1}\over{a_1(x)+}}\>{{1}\over{a_2(x)+}}\>{{1}\over{a_3(x)+}}\cdots\>.
\]

Let ${\cal B}$ denote Borel subsets of $(0,1]$ and ${\rm P}$ denote the Gauss' measure
\[ {\rm P}(A)={{1}\over{\ln{2}}}\int_A {{dx}\over{1+x}},\qquad A\in{\cal B}.\]
It is well-known that the random sequence $\{a_n\}$
defined on the probability space $((0,1], {\cal B}, {\rm P})$
is strictly stationary.

The literature concerning the limit theory
for functionals of digits of continued fraction
expansion (see e.g. \cite{Kn64}, \cite{IL71},
\cite{Io90}, \cite{IT90}, \cite{LL96}, \cite{Bi99})
reveals that some
classical results for sums of i.i.d. random variables
(see e.g. \cite{GnKo}, \cite{Lo60}, \cite{IL71}, \cite{Pet96})
found their full analogies in this theory.
One example is the Marcinkiewicz-Zygmund law of large numbers
(Cf. \cite{Ko33}, \cite{MZ37}),
or more generally the complete convergence
(Cf.  \cite{Sh95}, \cite{LL96}).

In this note we discover further analogies.
The first one is motivated by
Lemma in \cite{Fe45} (see also Theorem 4.13 in \cite{Pet96})
and \S4,  Ch. VI in \cite{Ko33}
(see also  Satz XII in \cite{Ko30} \& \cite{Ko28},
Theorem 1 in \S7, Ch. VII in \cite{Fe71}, Theorem 1.3 in \cite{Gu04}).

\begin{theorem}
\label{t0}
Let $c_n\to\infty$ be a sequence of positive numbers and $f$ be a Borel function.
In order that there exist a sequence $\{b_n\}$ such that
\begin{equation}
\label{ef}
c_n^{-1}(\sum_{k=1}^n f(a_k) - b_n )\rightarrow_{\rm P} 0
\end{equation}
it is necessary and sufficient that simultaneously
\begin{equation}
\label{em1}
n{\rm P}[\vert f(a_1)\vert > c_n]\to 0,\qquad
{{n}\over{c_n^2}}{\rm E}[f^2(a_1)I_{[\vert f(a_1)\vert\leq c_n]}]\to 0.
\end{equation}
If the latter conditions are satisfied, we can set
$b_n = n{\rm E}[f(a_1)I_{[\vert f(a_1)\vert\leq c_n]}]$.
\end{theorem}
Note that by the well-known formula
\begin{equation}
\label{ewkf}
x^2{\rm P}[\vert Z\vert>x] + {\rm E}[\vert Z\vert^2I_{[\vert Z\vert\leq x]}] =
 2\int_0^x y{\rm P}[\vert Z\vert>y]dy,
\end{equation}
(\ref{em1}) is equivalent to
\begin{equation}
\label{em}
\int_0^1 n{\rm P}[\vert f(a_1)\vert > \sqrt{x}c_n]dx\to 0.
\end{equation}
\vskip 2pt\noindent
If there exists some additional knowledge on
$c_n$ then (\ref{em1}) can be weakened.
\begin{theorem}
\label{t1}
Suppose that $c_n$ is a sequence of positive numbers such that
\begin{equation}
\label{eseq}
{\underline{\lim}}_{\>n\>}{{n}\over{c_{n}^2}}(c_{n+1}^2- c_{n}^2)>1.
\end{equation}
Then for any Borel function $f$
\begin{equation}
\label{ekf}
(c_n^{-1}\sum_{k=1}^n (f(a_k) - {\rm E}[f(a_k)I_{[\vert f(a_k)\vert\leq c_n]}]) \rightarrow_{\rm P} 0)
\Longleftrightarrow (n{\rm P}[\vert f(a_1)\vert>c_n]\rightarrow_n 0).
\end{equation}
\end{theorem}
Theorem \ref{t1} and Theorem 1.9.8 in \cite{BGT87} yield
\begin{corollary}
\label{c1}
Suppose that $f$ is a Borel function,
$c_n=n^{{1}\over{r}}h(n)$, where $h$ is a slowly
varying function in the sense of Karamata and $r\in(0,2)$.
Then the relation (\ref{ekf}) holds.
\end{corollary}
Let $c_n$ denote the
accumulated entrance fees up to the $n$-th trial
in the St. Petersburg game, i.e. $c_n=n\log_2{n}$ (Cf. \cite{Fe45}).
Since $n\ln{2}{\rm P}[a_1>n]\sim 1,$ therefore
${\rm E}[a_1I_{[a_1\leq c_n]}]\sim \log_2{n}$
and by Corollary \ref{c1} (see also Theorem 3 in \cite{Sz01})
we get for $f(x)=x$
\[ c_n^{-1}\sum_{k=1}^n a_k\to_{\rm P} 1. \]
On the other hand, by Theorem 4 in \cite{Sz01},
the convergence in probability cannot be replaced by the almost sure one.

The second analogy is motivated by the famous
characterization of the domain of attraction of the normal law
due to  Khinchine (Cf. \cite{Kh35}), L\'evy (Cf. \cite{Le35}) and Feller (Cf. \cite{Fe36})
(see also Corollary 1, \S5, Ch. XVII in \cite{Fe71})
\begin{theorem}
\label{t2}
Let $f$ be a Borel function.
In order that there exist sequences $\{c_n\}$ and $\{d_n\}$ such that
$\>\lim_{n\to\infty}{\cal L}(c_n^{-1}(\sum_{k=1}^n f(a_k)-d_n))= {\cal N}(0,1)$
it is necessary and sufficient that
the function ${\rm E}[f^2(a_1)I_{[\vert f(a_1)\vert\leq x]}]$
is slowly varying in the sense of Karamata. If the latter condition is satisfied, we can set
$d_n = n{\rm E}[f(a_1)]$.
\end{theorem}
Theorem 3 improves Proposition 2.9 (II) in \cite{Sa89}
obtained under additional assumptions on $c_n$ and $d_n$.
Furthermore, we no longer have to assume that condition $(*)$
on page 56 of \cite{Sa89} is satisfied.
By Theorem \ref{t2} and the
results in \cite{Sz01} we get that functionals of digits in
continued fraction expansion satisfy the Raikov principle (Cf.
Twierdzenie 4, \S28, Ch. V in \cite{GnKo}), namely
\begin{corollary}
\label{c2}
Suppose that $f$ is a Borel function and $c_n$ is a positive sequence. Then
$c_n^{-2}\sum_{k=1}^n f^2(a_k)\to_{\rm P} 1$ and ${\rm E}[f(a_1)]=0$
if and only if
\[ \>\lim_{n\to\infty}{\cal L}(c_n^{-1}\sum_{k=1}^n f(a_k))
= {\cal N}(0,1).\]
\end{corollary}
As it can be concluded from the rest of this note, the
presented results
remain true in a little bit more general mathematical
environment.

\section{Preliminaries}

We group here different results that will be used later on.
Let $\{X_k\}_{k\in{\mathbb N}}$ be a random sequence defined on a probability space $(\Omega, {\cal F}, P),$
denote by  $\{\tilde{X}_k\}$ it's independent copy
and in the case of stationarity by $\{X_k^{*}\}$ it's i.i.d. associated sequence
(all sequences are sharing the same probability space).
Define
\[ S_n=\sum_{k=1}^n X_k,\quad \widehat{X}_k = X_k - \tilde{X}_k,\quad
\widehat{S}_n = \sum_{k=1}^n \widehat{X}_k,\]
\[ M_n=\max_{1\leq k\leq n}\vert X_k\vert,\quad
\widehat{M}_n=\max_{1\leq k\leq n}\vert \widehat{X}_k\vert,\quad
M_n^{*}= \max_{1\leq k\leq n}\vert X_k^{*}\vert. \]
Denote by
${{\cal F}_{\!\!{k}}^{\,m}}$ the $\sigma$--field generated
$X_k, X_{k+1},\ldots, X_m$,
$m\in{\mathbb N}$, and recall the following
coefficients of dependence
\[\psi_n =
\sup_{k\in{\mathbb N}}
\sup\{
\vert {{P(A\cap B)}\over {P(A)P(B)}}-1\vert;\> P(A)P(B)>0,
\> A\in{{\cal F}_{\!\!{1}}^{\,k}}
,B\in {\cal F}_{n+k}^{\infty}\};
\]
\[\psi^{*}_n = \sup_{k\in{\mathbb N}}\sup\{
{{P(A\cap B)}\over {P(A)P(B)}};\> P(A)P(B)>0,
\> A\in{{\cal F}_{\!\!{1}}^{\,k}}
,B\in {\cal F}_{n+k}^{\infty}\}; \]
\[\psi'_n = \inf_{k\in{\mathbb N}}\inf\{
{{P(A\cap B)}\over {P(A)P(B)}};\> P(A)P(B)>0,
\> A\in{{\cal F}_{\!\!{1}}^{\,k}}
,B\in {\cal F}_{n+k}^{\infty}\}; \]
\[\varphi_n = \sup_{k\in{\mathbb N}}\sup\{
\vert P(B\,\vert\,A) - P(B) \vert;\> P(A)>0,
\> A\in{{\cal F}_{\!\!{1}}^{\,k}}
,B\in {\cal F}_{n+k}^{\infty}\}; \]
\[\rho_n =
\sup_{k\in{\mathbb N}}
\sup\{ \vert { Corr(f,g)}\vert;
\> f\in{\rm L}_{\rm real}^2({{\cal F}_{\!\!{1}}^{\,k}})\,
,g\in {\rm L}_{\rm real}^2({\cal F}_{n+k}^{\infty})\}.\]
It is well-known that
\begin{equation}
\label{edep}
\psi_n = \max\{\psi^{*}_n-1,1- \psi'_n\},\>
\psi_n\geq 2\varphi_n,\>
\psi_n\geq \rho_n,\>
1-\varphi_n\geq \psi'_n,\>
4\varphi_n\geq (\rho_n)^2,
\end{equation}
for every $n\geq 1$ (Cf. \cite{Br05}, p.109).
By Theorem 5.2 in \cite{Br05} the symmetrized $\psi_n,\psi^{*}_n,\psi'_n,\varphi_n,\rho_n$ coefficients of
the sequence $\{\widehat{X}_k\}$,
say $\widehat{\psi}_n, \widehat{\psi}^{*}_n, \widehat{\psi}'_n, \widehat{\varphi}_n, \widehat{\rho}_n,$
satisfy
\begin{equation}
\label{esd}
\widehat{\psi}_n\leq (1+\psi_n)^2-1,\>
\widehat{\psi}^{*}_n \leq (\psi^{*}_n)^2,\>
\widehat{\psi}'_n \geq (\psi'_n)^2,\>
\widehat{\varphi}_n \leq 1-(1-\varphi_n)^2,\>
\widehat{\rho}_n \leq \rho_n,
\end{equation}
for every $n\geq 1.$

The following lemma is a dependent version of L\'evy's inequality.
\begin{lemma}
\label{sum}
Assume that ${\cal L}(S_n-S_k)$ are symmetric for $n > k\geq 1$.
Then for $n\geq 1$ and $x>0$
\begin{eqnarray*}
2P[\vert S_n\vert > x]
\geq \psi'_1
P[\max_{1\leq k\leq n} \vert S_k\vert > x].
\end{eqnarray*}
\end{lemma}
\begin{proof} {\sc Lemma \ref{sum}}

Consider the sets
\[ C_1^{+} = [ S_1 > x ],\quad
C_k^{+} = [ \vert S_1\vert\leq x, \vert S_2\vert\leq x,\ldots ,
\vert S_{k-1}\vert\leq x, S_k > x ], \]
\[C_1^{-} = [ -S_1 > x ],\quad
   C_k^{-} = [\vert S_1\vert\leq x, \vert S_2\vert\leq x,\ldots ,
\vert S_{k-1}\vert\leq x, -S_k > x ], \]
and set $C_k=C_k^{+}\cup C_k^{-},$ $C=\bigcup_{k=1}^n C_k$.
Since the $C_k$ are disjoint we obtain
\begin{eqnarray*}
2P[\vert S_n\vert > x]
& = &
2P[\max_{1\leq k\leq n}\vert S_k\vert> x,\,
\vert S_n\vert > x]\\
& = &
2\sum_{k=1}^{n} P( C_k \cap [ \vert S_n \vert > x  ])   \\
& \geq & 2\sum_{k=1}^{n} (P( C_k^{+} \cap [ S_n > x   ]) +
P(C_k^{-} \cap [ -S_n > x ])) \\
& \geq &  2\sum_{k=1}^{n}
(P( C_k^{+} \cap [ S_n - S_k \geq 0 ])
+ P(C_k^{-} \cap [ -S_n + S_k \geq 0 ])) \\
& \geq &
\psi'_1
\sum_{k=1}^{n} (P(C_k^{+})
+
P(C_k^{-}))
=
\psi'_1
\sum_{k=1}^{n} P(C_k)\\
& = &
\psi'_1P(C)
=
\psi'_1
P[\max_{1\leq k\leq n} \vert S_k\vert > x],
\end{eqnarray*}
because $2P[S_n - S_k\geq 0] = 2P[S_n - S_k\leq 0] = 1 + P[S_n-S_k=0]\geq 1$.
This proves Lemma \ref{sum}.
\end{proof}

It is well-known that if ${ E}\vert X\vert<\infty,$ ${ E}\vert Y\vert<\infty,$
and $X$ is ${{\cal F}_{\!\!{1}}^{\,k}}\,$ measurable while
$Y$ is ${\cal F}_{n+k}^{\infty}$ measurable then
\begin{equation}
\label{epsi}
\vert{ E}[XY] - { E}[X]{ E}[Y]\vert \leq \psi_n{ E}\vert X\vert { E}\vert Y\vert
\end{equation}
(Cf. Lemma 1.2.11 in \cite{LL96}).

The statement below is a dependent version
of Kolmogorov's inverse inequality (Cf. \cite{Lo60}, p.235, \cite{AG80}, Theorem 2.8).
\begin{lemma}
\label{kol}
Let $\{X_k\}$ be a strictly stationary sequence such that
$\vert X_k\vert\leq c$  and ${\cal L}(S_k)$ are symmetric.
Then for every $x > 0$
\[ P[\max_{1\leq k\leq n}\vert S_k\vert> x]\geq
{{\psi'_1({ E}[S_n^2] - x^2)}\over{4(1+\psi_1){ E}[S_n^2] + \psi'_1(2(c+x)^2}-x^2)},\qquad n\geq 1.
\]
\end{lemma}
\begin{proof} {\sc Lemma \ref{kol}}

Let $C_k$ and $C$ be as in the proof of Lemma \ref{sum}.
We have
\[{ E}[S_n^2I_{C_k}]
\leq  2{ E}[(S_n-S_k)^2I_{C_k}]  + 2{ E}[S_k^2I_{C_k}]. \]
By (\ref{epsi}) we get
\[ { E}[(S_n-S_k)^2I_{C_k}] \leq (1+\psi_1){ E}[(S_n-S_k)^2]P(C_k).\]
Since $\vert S_kI_{C_k}\vert\leq (c+x)I_{C_k}$ so
by stationarity and Lemma \ref{sum} we obtain ($S_0=0$)
\begin{eqnarray*}
\lefteqn{
{ E}[S_n^2I_{C}]}\\
& &
\leq  2\sum_{k=1}^n((1+\psi_1){ E}[(S_n-S_k)^2]P(C_k)
+ (c+x)^2P(C_k))\\
& &
\leq  2\sum_{k=1}^nP(C_k)((1+\psi_1)\max_{1\leq k\leq n}{ E}[(S_n-S_k)^2]
+ (c+x)^2)\\
& &
=  2\sum_{k=1}^nP(C_k)((1+\psi_1)\max_{1\leq k< n}{ E}[S_k^2]
+ (c+x)^2)\\
& &
\leq  2\sum_{k=1}^nP(C_k)((1+\psi_1){ E}[\max_{1\leq k\leq n}S_k^2]
+ (c+x)^2)\\
& &
\leq  2\sum_{k=1}^nP(C_k)({{2(1+\psi_1)}\over{\psi'_1}}{ E}[S_n^2]
+ (c+x)^2)\\
& &
=  2P(C)({{2(1+\psi_1)}\over{\psi'_1}}{ E}[S_n^2] + (c+x)^2).
\end{eqnarray*}
On the other hand
\[ { E}[S_n^2I_{C}] = { E}[S_n^2]
- { E}[S_n^2I_{\Omega\setminus C}]\geq { E}[S_n^2] + x^2P(C) -x^2.\]
This completes the proof.
\end{proof}

The next inequality follows from the proof on p.298 in \cite{Pe90}.
\begin{proposition}
\label{ap3}
Suppose $\{X_k\}$ is a strictly stationary sequence and $\varphi_m<1$.
Then for every $x\geq 0$ and every $n\geq m\geq 1$
\[
(1-\varphi_m)P[M_{\lfloor n/m \rfloor}^{*}>x]\leq P[M_n>x]
\leq m(1+\varphi_m) P[M_{\lfloor n/m \rfloor+1}^{*}>x].
\]
\end{proposition}

We will need the following estimate
\begin{lemma}
\label{cut}
Suppose $\{X_k\}$ is a strictly stationary sequence and
$\psi'_1>0$.
Then for every $x\geq 0$ and every  $n\geq m\geq 1$
\[(1-\widehat{\varphi}_1)\widehat{\psi}'_1
(1 - e^{-nP[\widehat{X}_1 > 2x]})
\leq 4m
P[\vert\widehat{S}_{\lfloor n/m \rfloor+1}\vert > x]. \]
\end{lemma}
\begin{proof} {\sc Lemma \ref{cut}}

It is easy to see that for every $x\geq 0$ we have
\begin{equation}
\label{ems}
P[M_n > 2x ] \leq P[\max_{1\leq k\leq n}\vert S_k\vert > x].
\end{equation}
Therefore by Proposition \ref{ap3} and Lemma \ref{sum} we obtain
\begin{eqnarray*}
\lefteqn{
(1-\widehat{\varphi}_1)\widehat{\psi}'_1(1 - e^{-nP[\widehat{X}_1 > 2x]})\leq
(1-\widehat{\varphi}_1)\widehat{\psi}'_1P[\max_{1\leq k\leq n}\vert\widehat{X}^{*}_k\vert > 2x]}\nonumber\\
& & \leq
\widehat{\psi}'_1P[\widehat{M}_n > 2x]
\leq 2m\widehat{\psi}'_1P[\widehat{M}_{\lfloor n/m \rfloor+1}>2x]\nonumber\\
& &
\leq 2m\widehat{\psi}'_1P[\max_{1\leq k\leq \lfloor n/m \rfloor+1}\vert \widehat{S}_k\vert > x]
\leq 4m
P[\vert\widehat{S}_{\lfloor n/m \rfloor+1}\vert > x].
\end{eqnarray*}
This is our assertion.
\end{proof}

The following estimates are consequences of Lemma 2.1 and 2.3 in \cite{Br88}.
\begin{proposition}
\label{var}
Suppose $\rho_1<1$ and $\sum_{n=1}^{\infty}\rho_{2^n}<\infty$ for a strictly stationary sequence
$\{X_k\}$ with finite variances.
Then there exist positive constants $C,D$ depending only on $\{\rho_k\}$ such
that $\>Cn{ Var}[X_1]\leq {Var}[S_n]\leq Dn{ Var}[X_1],$ $n\geq 1$.
\end{proposition}

The proof of the next statement is an easy consequence of the hint on p.91 in \cite{AG80} and
it is included here for the reader convenience.
\begin{lemma}
\label{mom}
Suppose that $\{X_k\}$ is a strictly stationary sequence such that $\psi'_1>0$ and
${\cal L}(c_n^{-1}(S_n - b_n))$
are asymptotically normal for
$c_n^2=nh(n)$
where $h(n)$ is a slowly varying sequence. Then
$\sup_n\>c_n^{-q}{ E}\vert \widehat{S}_n\vert^q < \infty,$
for any $q\in(0,2)$.
\end{lemma}
\begin{proof} {\sc Lemma \ref{mom}}

Suppose $d>0$ and ${\delta}
\in(0,{{1}\over{2}}\widehat{\psi}'_1(1-{\widehat{\varphi}}_1))$
are such that
\[ P[c_n^{-1}\vert \widehat{S}_n\vert > d] \leq {\delta}, \]
for $n>N_{\delta}$.
By Lemma \ref{sum} we have for $n > N_{\delta}$
\[ P [c_{mn}^{-1}\max_{1\leq k\leq m}
\vert\widehat{S}_{nk} - \widehat{S}_{n(k-1)}\vert
> d ] \leq
{{2}\over{\widehat{\psi}'_1}}
P[c_{mn}^{-1}\vert\widehat{S}_{mn}\vert > d]
\leq
{{2\delta}\over {\widehat{\psi}'_1}}\cdotp \]
On the other hand by Proposition \ref{ap3}
\begin{eqnarray*}
\lefteqn{
P[c_{mn}^{-1}\max_{1\leq k\leq m}
\vert\widehat{S}_{nk} - \widehat{S}_{n(k-1)}\vert
> d]}\\
& & \geq
(1-{\widehat{\varphi}}_1) P[c_{mn}^{-1}
\max_{1\leq k\leq m}
\vert\widehat{S}_{n,k}^{*}\vert > d)
\\&  &
=   (1-{\widehat{\varphi}}_1)(1 - (1 - P[
c_{mn}^{-1}\vert\widehat{S}_n\vert > d]
)^m),
\end{eqnarray*}
where $\widehat{S}_{n,k}^{*}$ are independent copies of
$\widehat{S}_{nk} - \widehat{S}_{n(k-1)}$. Thus
\[ P[c_{mn}^{-1}\vert\widehat{S}_n\vert > d]\leq
1 - (1 - {{2\delta}\over{\widehat{\psi}'_1(1-{\widehat{\varphi}}_1)}})^{1\over m}. \]
From
\[ 1-(1-a)^{1\over m}\sim {1\over m}\ln{1\over {1-a}},\quad
m\rightarrow \infty,\quad  0< a<1, \]
it follows that there exists a constant
$A=A({\delta}, \widehat{\psi}'_1, \widehat{\varphi}_1, N_{\delta})$ such that
\begin{equation}
\label{je12}
P[c_{mn}^{-1}\vert\widehat{S}_n\vert > d] \leq
{A\over m}\cdotp
\end{equation}
Since $c_n^2=nh(n),$  where $h(n)$ is a slowly varying sequence
so by the Uniform Convergence Theorem for R
(Cf.  \cite{BGT87}, Theorem 1.5.2, p.22)
we have that for every
${\gamma}>0$ there exists $N_{\gamma} >N_{\delta}$ such that
\[ \bigg\vert {{m^{-\gamma}n^{-\gamma}h(mn)}\over{n^{-\gamma}h(n)}} -
 m^{-\gamma} \bigg\vert  = \bigg\vert {{c_{mn}}\over{m^{{1\over 2}
+{\gamma}}c_n}} - {1\over {m^{\gamma}}}\bigg\vert < {1\over 2}, \]
for $n\geq N_{\gamma}$ uniformly in $m \geq 1$.
Thus
\[ {{c_{mn}}\over{c_n}}  < m^{{1\over 2}+{\gamma}}(m^{-\gamma} +
{1\over 2}) \leq {3\over 2} m^{{1\over 2}+{\gamma}}, \]
for $n\geq N_{\gamma}$,  $m\geq 1$ which with (\ref{je12}) gives
\[ m P[c_n^{-1}\vert\widehat{S}_n\vert
> {3\over 2}dm^{{1\over 2}+{\gamma}}]
\leq A, \]
for $n\geq N_{\gamma}$, $m\geq 1$. Replacing $m$ with $m^{2\over{1+2{\gamma}}}$
we obtain
\[  m^{2\over{1+2{\gamma}}} P[c_n^{-1}\vert\widehat{S}_n\vert > md']
\leq {A'}, \]
where ${A'}= {A'}(A,{\gamma}),$  ${d'}={d'}(d, {\gamma})$ and $n\geq N_{\gamma}$,
$m\geq 1$. Now, taking ${\gamma}$ such that
$ 0< {\gamma} <  {1\over q} - {1\over 2}$ we get
\[ m^{-1+q} P[c_n^{-1}\vert\widehat{S}_n\vert > m{d'}] \leq
{{A'}\over{m^{1+c}}}, \]
where $c={2\over {1+2{\gamma}}} - q >0$.
Thus
\[ \sum_{m=1}^{\infty} m^{q-1}
P[ c_n^{-1}\vert\widehat{S}_n\vert > m{d'}] < \infty, \]
for $n\geq N_{\gamma}$.
This completes the proof.
\end{proof}

The next inequality follows from Proposition \ref{ap3} and Proposition 4.3, Lemma 4.2
in \cite{HK98}.
\begin{proposition}
\label{mmom}
Suppose that $\{X_k\}$ is a strictly stationary random sequence such that
$\phi_1<1$ and $\sup\{x; P[\vert X_1\vert \leq x]<1\}=\infty$. Then for every $q,x>0$
\[(1-\phi_1){{nE[\vert X_1\vert^qI_{[\vert X_1\vert>x]}]}\over{1+nP[\vert X_1\vert >x]}}
\leq E[M_n^q],\qquad n\geq 1. \]
\end{proposition}

The lemma below is a dependent analog of Khinchine's inequality (Cf. \cite{AG80}, p.176).
\begin{lemma}
\label{khin}
Suppose that $\{X_k\}$ is a random sequence such that $\vert X_k\vert<c$ and
$0<\psi'_1\leq \psi^{*}_1<\infty$. Then for any $q\in[1,p)$ there exists a constant
$B_{pq}$ depending only on $\psi'_1,\psi^{*}_1,p,q$ such that
\[ { E}\vert \widehat{S}_n\vert^p \leq B_{pq}
\max\{{ E}^{p\over q}\vert  \widehat{S}_n\vert^q, c^p\},\qquad n\geq 1.\]
\end{lemma}
\begin{proof} {\sc Lemma \ref{khin}}

From the Nagaev generalization of the inequality (3.3) in \cite{HJ74}
(Cf. the relation $(11)$
and the proof of Lemma in \cite{Na00}) it follows that for
any random
sequence $\{ X_k\}$ such that $\vert X_k\vert<c$
we have for $t\geq c$
\begin{equation}
\label{enmax}
P[\max_{1\leq k\leq n}\vert S_k\vert>4t]\leq \psi^{*}_1
(P[\max_{1\leq k\leq n}\vert S_k\vert>t])^2.
\end{equation}
In view of this and Lemma \ref{sum} we get for $t>2c$
\[ P[\vert \widehat{S}_n\vert > 4t]\leq (
\gamma P[\vert \widehat{S}_n\vert >t])^2, \]
where $\gamma = {{2\sqrt{\widehat{\psi}^{*}_1}}\over{\widehat{\psi}'_1}}\cdotp$
If we set
\[ t_0 = 4^{1\over q}{\gamma}^2\max\{ { E}^{1\over q}\vert \widehat{S}_n\vert^q, c\} > 2c \]
then by the Markov inequality we obtain for $t\geq t_0$
\[
{\gamma}^2P[\vert \widehat{S}_n\vert >t] \leq {1\over 4}\cdotp \]
From this we have
\begin{eqnarray*}
{ E}\vert \widehat{S}_n\vert^p & = & p\int_0^{4t_0} x^{p-1}
P[\vert \widehat{S}_n\vert > x] dx + \\
& & \quad + p\sum_{k=1}^{\infty}\int_{{4^k}t_0}^{{4^{k+1}}t_0} x^{p-1}
P[\vert \widehat{S}_n\vert > x] dx \\
& \leq & (4t_0)^p + t_0^p\sum_{k=1}^{\infty}
4^{p(k+1)}P[\vert \widehat{S}_n\vert > 4^k t_0] \\
& \leq & (4t_0)^p +  4^pt_0^p
\sum_{k=1}^{\infty}  4^{kp}{\gamma}^{\sum_{i=1}^k 2^i}
P^{2^k}[\vert \widehat{S}_n\vert > t_0] \\
& \leq & (4t_0)^p +
 4^pt_0^p{\gamma}^{-2}
\sum_{k=1}^{\infty}
4^{kp}({\gamma}^2P[\vert \widehat{S}_n\vert > t_0])^{2^k} \\
& \leq & (4t_0)^p +
4^pt_0^p{\gamma}^{-2}
\sum_{k=1}^{\infty} {{4^{kp}}\over{4^{2^k}}} \\
& = & (4t_0)^p(1 + {1\over{{\gamma}^2}})
\sum_{k=1}^{\infty} {{4^{kp}}\over{4^{2^k}}} \\
& \leq & 4^{p\over  q}2{\gamma}^{2p}4^p
\sum_{k=1}^{\infty} {{4^{kp}}\over{4^{2^k}}}
\max\{{ E}^{p\over q}\vert \widehat{S}_n\vert^q, c^p \} \\
& = & B_{pq}
\max\{{ E}^{p\over q}\vert  \widehat{S}_n\vert^q, c^p\}.
\end{eqnarray*}
This is desired conclusion.
\end{proof}

We will also need the following symmetrization result for slowly varying functions.
\begin{lemma}
\label{sym}
If ${ E}[\widehat{X}^2 I_{[\vert\widehat{X}\vert\leq x]}]$ varies slowly
then ${ E}[X^2 I_{[\vert X\vert\leq x]}]$ varies slowly, too.
\end{lemma}
\begin{proof} {\sc Lemma \ref{sym}}

Since ${ E}\vert\widehat X\vert<\infty$ thus ${ E}\vert X\vert<\infty$
(Cf. \cite{Lo60}, p.243).
By $P[\vert X\vert \geq 2{ E}\vert X\vert]
\leq {{1}\over{2}}$ we have
that $median(X)\leq 2E\vert X\vert$.
Hence by the weak symmetrization
inequalities (Cf. \cite{Lo60}, p.245) we have
\[
{{1}\over{2}} P[\vert X\vert > x + 2{ E}\vert X\vert]
\leq
P[\vert\widehat X\vert > x]
\leq
2P[\vert X \vert > {{1}\over{2}}x].
\]
Therefore, if  $x\geq e{ E}\vert X\vert$ then
\begin{eqnarray}
\label{esym}
{{x^2P[\vert X\vert> x]}\over{2\int_0^x yP[\vert X\vert> y]dy}}
& \leq &
{{16x^2}\over{(x-2{ E}\vert X\vert)^2}}
{{(x-2{ E}\vert X\vert)^2 P[\vert \widehat X\vert
> x-2{E}\vert X\vert]}
\over{2\int_0^{x-2E\vert X\vert}
yP[\vert \widehat X\vert > y]dy}} \nonumber\\
& & \qquad\qquad\qquad \times
{{\int_0^{x-2E\vert X\vert}
yP[\vert \widehat X\vert > y]dy}
\over{\int_0^{2x}
yP[\vert \widehat X\vert > y]dy}}\cdotp
\end{eqnarray}
Since the fraction standing left of the formula number (\ref{esym}) is at most 1
thus by (\ref{ewkf}),(\ref{esym}) and Theorem 2, VIII, \S9 in \cite{Fe71}
\[
{{x^2P[\vert X\vert > x]}\over{x^2P[\vert X\vert > x]
+ E[X_1^2I_{[\vert X\vert\leq x]}]}} =
 {{x^2P[\vert X\vert > x]}\over{2\int_0^x yP[\vert X\vert> y]dy}}
\to 0.
\]
This proves Lemma \ref{sym}.
\end{proof}

By Lemma 2.1 in \cite{Ph88} and Corollary in \cite{Io89}
the sequence $\{ f(a_k) \}$ fulfills
$\psi_n\leq \varrho^n$,
for some
$\varrho<0.8$ and $\psi_1\leq 2\ln{2}-1 < 0.39$.
Therefore, by (\ref{edep}) and (\ref{esd}) we have
$\varphi_1< 1/2$,
$\widehat{\varphi}_1 < 1/2,$ ${\psi}'_1 > 0,$
$\widehat{\psi}'_1 > 0,$
$\psi^{*}_1 < 2,$
$\widehat{\psi}^{*}_1 < 2,$
$\widehat{\rho}_1\leq \widehat{\psi}_1< 1$.

In particular, by Lemma \ref{sum} and (\ref{ems})
we have for every $x\geq 0$ and $n\geq 1$
\begin{equation}
\label{emax}
{\rm P}[\max_{1\leq k\leq n}\vert \widehat{f}(a_k)\vert > 2x]
\leq {\rm P}[\max_{1\leq k\leq n}\vert \widehat{S}_k\vert > x]
\leq 6{\rm P}[\vert \widehat{S}_n\vert > x].
\end{equation}

\section{Proofs}
\vskip 2pt\noindent
{\sc Proof of Theorem \ref{t0} }

Without the loss of generality we may assume that
${{\rm E}[f^2(a_1)]}=\infty$ (the case ${{\rm E}[f^2(a_1)]}<\infty$ is covered by Theorem \ref{t2}).
Assume first that (\ref{em1}) holds.
Set $X_{nk} = f(a_k)I_{[\vert f(a_k)\vert\leq c_n]}$.
By Proposition \ref{var} we obtain
\[ {\rm E}(\sum_{k=1}^n (X_{nk} - {\rm E}[X_{n1}]))^2\!
\leq Dn{\rm Var}[X_{n1}]. \]
In view of this and Chebyshev's inequality we get for any $\epsilon>0$
\begin{eqnarray*}
\lefteqn{
{\rm P}[\vert \sum_{k=1}^n (f(a_k) -{\rm E}[f(a_k)I_{[\vert f(a_k)\vert\leq c_n]}] )\vert > \epsilon c_n ]}\\
& &\leq
{\rm P}[\vert\sum_{k=1}^nf(a_k)I_{[\vert f(a_k)\vert > c_n ]}\vert>{{\epsilon}\over{2}}c_n]
+ {\rm P}[\vert\sum_{k=1}^n (X_{nk} - {\rm E}[X_{nk}])\vert >{{\epsilon}\over{2}}c_n]\\
& &
\leq n{\rm P}[\vert f(a_1)\vert > c_n]
 +
8Dn\epsilon^{-2}c_n^{-2}{\rm E}[f(a_1)^2I_{[\vert f(a_1)\vert\leq c_n]}]\to 0.
\end{eqnarray*}

Conversely, assume (\ref{ef}).
By (\ref{emax})
\[ {\rm P}[\max_{1\leq k\leq n}\vert\widehat{f}(a_k)\vert > \epsilon c_n]\to 0  \]
so that denoting
${Y}_{nk}(\epsilon) = \widehat{f}(a_k)I_{[\vert\widehat{f}(a_k)\vert\leq \epsilon c_n]}$
we get
\[ c_n^{-1}\sum_{k=1}^n Y_{nk}(\epsilon)\to_{\rm P} 0,\qquad \epsilon\in(0,1).\]
Let $Z_n(\epsilon)=\sum_{k=1}^n Y_{nk}(\epsilon)$
then $c_{n}^{-1}\widehat{Z}_n(\epsilon)\to_{\rm P} 0$.
Thus by Lemma \ref{sum}
\[ c_n^{-1}\max_{1\leq k\leq n}\sum_{m=1}^k \widehat{Y}_{nm}(\epsilon)\to_{\rm P} 0\]
and by Lemma \ref{kol} (with $x=\epsilon, X_k = c_n^{-1}\widehat{Y}_{nk}$)
we get
\[
\lim c_n^{-2}{\rm E}[\widehat{Z}_n^2(\epsilon)]=0.\]
Since ${\rm E}[\widehat{f}^2(a_1)]=2{\rm E}[f^2(a_1)]=\infty$
thus ${\rm E}^2\vert Y_{n1}(\epsilon)\vert  = o({\rm E}[Y_{n1}^2(\epsilon)])$
(Cf. (2.6.14) in \cite{IL71}).
Therefore by Proposition \ref{var}
\[ \lim_n {{n}\over{c_n^2}}{\rm E}[Y_{n1}^2(\epsilon)]
=\lim_n {{n}\over{c_n^2}}{\rm Var}[Y_{n1}(\epsilon)]
=\lim_n {{n}\over{2c_n^2}}{\rm Var}[\widehat{Y}_{n1}(\epsilon)]=0,\quad \epsilon\in(0,1) \]
so that $n^{-1}c_n^2\to\infty$.
Further, by (\ref{ewkf}) and since
$n{\rm P}[\vert \widehat{f}(a_1)\vert > \epsilon c_n]\to 0$ hence
\[ \int_0^1 n{\rm P}[\vert f(a_1) - \tilde{f}(a_1)\vert > \sqrt{x\epsilon} c_n] dx \to 0,
\qquad \epsilon\in(0,1). \]
Now, by the weak symmetrization
inequalities (Cf. \cite{Lo60}, p.245) we have
\begin{eqnarray*}
\lefteqn{
\int_0^1 n{\rm P}[\vert f(a_1)\vert > \sqrt{x}c_n] dx }\\
& &
\leq \int_0^1 n{\rm P}[\vert f(a_1) - median(f(a_1))\vert > \sqrt{x}{{c_n}\over{2}}] dx\\
& &
\qquad + \int_0^1 n{\rm P}[\vert median(f(a_1))\vert > \sqrt{x}{{c_n}\over{2}}] dx\\
& &
\leq \int_0^1 2n{\rm P}[\vert f(a_1) - \tilde{f}(a_1)\vert > \sqrt{x}{{c_n}\over{2}}] dx \\
& & \qquad
+ \int_0^1 n{\rm P}[\vert median(f(a_1))\vert > \sqrt{x}{{c_n}\over{2}}] dx\\
& &
\leq \int_0^1 2n{\rm P}[\vert f(a_1) -  \tilde{f}(a_1)\vert > \sqrt{x}{{c_n}\over{2}}] dx
+ 4{{n}\over{c_n^2}}(median(f(a_1)))^2\to 0.
\end{eqnarray*}
This proves Theorem \ref{t0}.

\vskip 5pt\noindent
{\sc Proof of Theorem \ref{t1} }

By (\ref{eseq}) there exists $N\in{\mathbb N}$ such that for $n>N$
\[ {{n-1}\over{c_{n-1}^2}}(c_n^2- c_{n-1}^2)\geq c >1. \]
Setting, if necessary $c_k^2=(c+1)^{-N+k+1}c_{N-1}^2$
for $1\leq k < N,$ we may assume that
\begin{equation}
\label{eseq1}
{{n-1}\over{c_{n-1}^2}}(c_n^2- c_{n-1}^2)\geq c,\quad n>1,\quad c>1,
\end{equation}
or, equivalently
\begin{equation}
\label{eseq2}
{{c_n^2}\over{n}} - {{c_{n-1}^2}\over{n-1}}
\geq {{c-1}\over{c}}{{c_n^2-c_{n-1}^2}\over{n}},\quad n>1,\quad c>1.
\end{equation}
Further, by (\ref{eseq1}) we have
\[\ln{c_n^2}\geq  \ln{c_1^2} + \sum_{k=2}^n (\ln(1+{{c}\over{k-1}})-{{c}\over{k-1}})
+ c\sum_{k=2}^n {{1}\over{k-1}}\]
and since by the Taylor expansion of $\ln{(1+x)}$ the first sum on the right hand side is $O(1)$
so that ${\underline{\lim}}_{\>n} {{\ln{{{c_n^2}\over{n}}}}\over{\ln{n}}}\geq c-1>0$.
By this, (\ref{eseq1}) and  (\ref{eseq2}) we get  ${{c_n^2}\over{n}}\nearrow \infty$.

Using in (\ref{eseq2}) the convention ${{0}\over{0}}=0$ we have
\begin{eqnarray*}
\lefteqn{
{{n}\over{c_n^2}}{\rm E}[f(a_1)^2I_{[\vert f(a_1)\vert\leq c_n]}]}\\
& & \leq {{n}\over{c_n^2}}
\sum\limits_{k=1}^n c_k^2 {\rm P}[\vert f(a_1)\vert\in(c_{k-1},c_k]]\\
& & = {{n}\over{c_n^2}}
\sum\limits_{k=1}^n\sum\limits_{\nu=1}^k (c_{\nu}^2-c_{\nu-1}^2) {\rm P}[\vert f(a_1)\vert\in(c_{k-1},c_k]]\\
& & = {{n}\over{c_n^2}}
\sum\limits_{\nu=1}^n (c_{\nu}^2-c_{\nu-1}^2)\sum\limits_{k=\nu}^n {\rm P}[\vert f(a_1)\vert\in(c_{k-1},c_k]]\\
& & \leq  {{n}\over{c_n^2}}
\sum\limits_{\nu=1}^n ({{c_{\nu}^2-c_{\nu-1}^2}\over{\nu}}) \nu {\rm P}[\vert f(a_1)\vert > c_{\nu-1}]\\
& & \leq
{{c}\over{c-1}}{{n}\over{c_n^2}}
\sum\limits_{\nu=1}^n ({{c_{\nu}^2}\over{\nu}}-{{c_{\nu-1}^2}\over{\nu-1}} )
\nu {\rm P}[\vert f(a_1)\vert > c_{\nu-1}]
\to_n 0
\end{eqnarray*}
by the Toeplitz lemma (Cf. \cite{Lo60}, p.238).
In view of Theorem \ref{t0}
the relation (\ref{ekf}) is proved.
\vskip 2pt\noindent
{\sc Proof of Theorem \ref{t2} }

Suppose that ${\rm E}[f^2(a_1)I_{[\vert f(a_1)\vert\leq x]}]$
is slowly varying and ${\rm E}[f(a_1)]=0$.
Define the sequence $b_n$ as follows:
{ if $\>{\rm E}[f^2(a_1)]<\infty$ then
$b_n = \sigma\sqrt{n}$ provided that
\[ \infty > \sigma^2 = {\rm E}[f^2(a_1)]
+ 2\sum_{k=2}^{\infty} {\rm E}[f(a_1)f(a_k)]>0; \]
if $\>{\rm E}[f^2(a_1)]=\infty$ then
\[ b_n = \sup\{x>0\, ;\, x^{-2}{\rm E}[\vert f(a_1)\vert^2I_{[\vert f(a_1)\vert\leq x]}]
\geq {{1}\over{n}} \} . \]}
Assume first that ${\rm E}[f^2(a_1)]<\infty$. Thus by
Theorem 18.5.2 in \cite{IL71}
and Proposition \ref{var}
we have $0<\sigma^2<\infty$
and
\[ {\cal L}(b_n^{-1}\sum_{k=1}^n f(a_k))\to {\cal N}(0,1).\]
Now, assume ${\rm E}[f^2(a_1)]=\infty$.
Then
$b_n^2\sim n{\rm E}[f^2(a_1)I_{[\vert f(a_1)\vert\leq b_n]}]$
and $b_n\to\infty$.
Let
$X_{nk}=f(a_k)I_{[\vert f(a_k)\vert \leq b_n]},$
$S_n=\sum_{k=1}^n X_{nk}$.
By the results in \cite{Sz92} if the following conditions are satisfied
\begin{equation}
\label{we1}
\lim_{n\to \infty}n {\rm P}[\vert f(a_n) \vert >{b_n}] = 0,
\end{equation}
\begin{equation}
\label{we2}
{\tau}_n^2 = {\rm Var}(\sum_{k=1}^n X_{nk})
\to \infty,
\end{equation}
\begin{equation}
\label{we3}
\lim_{n\to \infty}
{{\tau}_n^{-2}} {\rm E}[\max_{1\leq k\leq n}(X_{nk} - {\rm E}[X_{nk}])^2] = 0
\end{equation}
then
\begin{equation}
\label{we4} {\cal L}({\tau}_n^{-1}(S_n- n{\rm E}[X_{n1}])) \to
{\cal N}(0,1).
\end{equation}
The condition (\ref{we1}) easily follows by the definition of
$b_n$ and the slow variation of
${\rm E}[f^2(a_1)I_{[\vert f(a_1)\vert\leq x]}]$
(Cf. Theorem 2, VIII, \S9 in \cite{Fe71}).
For (\ref{we2}) let us observe that by (\ref{epsi}) and
since ${\rm E}[f^2(a_1)]=\infty$ thus
${\rm E}^2\vert X_{n1}\vert  = o({\rm E}[X_{n1}^2])$
(Cf. (2.6.14) in \cite{IL71}) and therefore
 we have $b_n^2\sim n{\rm Var}[X_{n1}]$. Now by
\begin{eqnarray*}
\lefteqn{
\tau_n^2 = {\rm E}(\sum_{k=1}^n (X_{nk} - {\rm E}[X_{n1}]))^2}\\
& & = n{\rm Var}[X_{n1}] + 2\sum_{k=2}^{n}(n-k+1) {\rm Cov}[X_{n1}X_{nk}]\\
& & = n{\rm Var}[X_{n1}](1 + O({{2}\over{(1-\varrho)}}{{n}\over{b_n^2}}{\rm E}^2\vert X_{n1}\vert))
= b_n^2(1+o(1))
\end{eqnarray*}
(\ref{we2}) follows.
For the relation (\ref{we3}) note that for every $\epsilon\in(0,1),$
some $K>0$ independent of $n$ by
the slow variation of
${\rm E}[f^2(a_1)I_{[\vert f(a_1)\vert\leq x]}]$
\begin{eqnarray*}
\lefteqn{
{\tau}_n^{-2} {\rm E}[\max_{1\leq k\leq n}(X_{nk} - {\rm E}[X_{nk}])^2]}\\
& &
\leq K(b_n^{-2}n{\rm E}[f^2(a_1)I_{[\vert f(a_1)\vert\leq b_n]}
I_{[f^2(a_1)I_{[\vert f(a_1)\vert\leq b_n]} >\epsilon b_n^2]}] + \epsilon)\\
& &
\leq K(b_n^{-2}n{\rm E}[f^2(a_1)I_{[\vert f(a_1)\vert\in (\sqrt{\epsilon}b_n, b_n]}] + \epsilon)
\leq K(o(1)+\epsilon).
\end{eqnarray*}
Thus (\ref{we4}) holds and since ${\rm E}[f(a_1)]=0,$ $\tau_n\sim b_n$ so that
${\cal L}(b_n^{-1}S_n)\to {\cal N}(0,1)$.
If ${\rm E}[f(a_1)]\neq 0$ then we consider the sequence $\{f(a_k)-{\rm E}[f(a_1)]\}$.
It is worth noting that by Theorem 1.9.8 in \cite{BGT87}
the sequence $\{b_n\}$ can be replaced by $\{c_n\}$ such that  $\lim_n {{n}\over{c_n^2}}(c_{n+1}^2-c_n^2)=1$.

Conversely, assume
${\phi}_n(c_n^{-1}\theta)e^{-i\theta{{d_n}\over{c_n}}}={\rm E}[e^{i\theta c_n^{-1}\sum_{k=1}^n{f}(a_k)}]
e^{-i\theta{{d_n}\over{c_n}}}\to e^{-{{\theta^2}\over{2}}}$
thus we have
$\widehat{\phi}_n(c_n^{-1}\theta)={\rm E}[e^{i\theta c_n^{-1}\sum_{k=1}^n\widehat{f}(a_k)}]\to e^{-\theta^2}$.
By Theorem 3.1 in \cite{Ja93} we have $c_n^2=nh(n)$, where
$h(n)$ is a slowly varying sequence.
There is no loss of generality in assuming that
${\rm E}[f^2(a_1)]=\infty$ so that $\sup\{x; {\rm P}[\vert \widehat{f}(a_1)\vert \leq x]<1\}=\infty$.
By Lemma \ref{mom} we have $\sup_n c_n^{-q}{\rm E}\vert\sum_{k=1}^n \widehat{f}(a_k)\vert^q<\infty,$
$q\in(0,2),$ so
by (\ref{emax})
we obtain that the sequence $\{c_n^{-1}\max_{1\leq k\leq n}\vert \widehat{f}(a_k)\vert\}$
is uniformly integrable.
On the other hand from Lemma \ref{cut} we get that
\begin{eqnarray*}
\lefteqn{
(1-\widehat{\varphi}_1)\widehat{\psi}'_1
\lim_n (1 - e^{-n{\rm P}[\widehat{f}(a_1) > \epsilon c_n]})}\\
& & \leq 4m\lim_n
{\rm P}[\vert\widehat{S}_{\lfloor n/m \rfloor+1}\vert > {{\epsilon}\over{2}}c_n]
= {{4m}\over{\sqrt{\pi}}}\int_{{{\epsilon}\over{2}}\sqrt{m}}^{\infty}  e^{-u^2/4}du
\end{eqnarray*}
and letting $m\to\infty$ it yields $\{c_n^{-1}\max_{1\leq k\leq n}\vert \widehat{f}(a_k)\vert\}\to 0$
in probability
since the integral tends to $0$ faster than exponentially.
Moreover, by Theorem 3.5 in \cite{Bi99} the latter convergence takes place in $L^1$, too.
Now, by Proposition \ref{mmom}  (with $q=1, x=\epsilon, X_k = c_n^{-1}\widehat{f}(a_k)$)
we obtain
\[ {{n}\over{c_n}}{\rm E}[\vert\widehat{f}(a_1)\vert I_{[\vert \widehat{f}(a_1)\vert > \epsilon c_n]}]
\to 0.\]
Hence
using the notation from the proof of Theorem \ref{t0} we get
\begin{eqnarray*}
\lefteqn{
\sup_n c_n^{-1}{\rm E}\vert Z_n(\epsilon)\vert}\\
& & \leq  \sup_n c_n^{-1}{\rm E}\vert \sum_{k=1}^n \widehat{f}(a_k)\vert
+ \sup_n c_n^{-1}{\rm E}\vert \sum_{k=1}^n \widehat {f}(a_k)I_{[\vert\widehat{f}(a_k)\vert>\epsilon c_n]}\vert\\
& & \leq  \sup_n c_n^{-1}{\rm E}\vert \sum_{k=1}^n \widehat{f}(a_k)\vert
+ \sup_n {{n}\over{c_n}}{\rm E}[\vert\widehat {f}(a_1)\vert I_{[\vert\widehat{f}(a_1)\vert>\epsilon c_n]}]<\infty.
\end{eqnarray*}
Consequently
$\sup_n c_n^{-1}{\rm E}\vert \widehat{Z}_n(\epsilon)\vert<\infty$
and Lemma \ref{khin} (with $p=2, q=1$) gives
$\sup_n c_n^{-2}{\rm E}[\widehat{Z}_n^2(\epsilon)]<\infty$.
In view of this it follows from Proposition \ref{var} that
$\sup_n {{n}\over{c_n^2}}{\rm E}[Y_{n1}^2(\epsilon)]<\infty$.
Therefore
there have to be $n^{-1}c_n^2\to\infty$
since we assumed
${\rm E}[f^2(a_1)]=\infty$.

By Lemma 1 in \cite{He87}  for $\theta\in{\Bbb R}$
and any integers $m>0,$ $p>1$ such that
\begin{equation}
\label{eh1}
{\rm E}\vert e^{i\theta c_n^{-1}\widehat{f}(a_1)} - 1\vert \leq
\min\{{{1}\over{2(1+\widehat{\psi}_1)^2(2m+1)^2}}, {{1}\over{2(1+\widehat{\psi}_1)(2pm+1)}} \}
\end{equation}
and
\begin{equation}
\label{eh2}
n{\rm E}^2\vert e^{i\theta c_n^{-1}\widehat{f}(a_1)} - 1\vert
\leq {{1}\over{2}}(9+\sum_{\nu=1}^m\widehat{\psi}_{\nu})^{-1}
n(1-\widehat{\phi}_1(c_n^{-1}\theta))
\end{equation}
we have
\begin{eqnarray}
\label{ehe}
\lefteqn{
\vert\widehat{\phi}_n(c_n^{-1}\theta) - \exp\{n(\widehat{\phi}_1(c_n^{-1}\theta) -1)\}\vert}\\
& & \leq
(9 + \sum_{\nu=1}^m\widehat{\psi}_{\nu})
n{\rm E}^2\vert e^{i\theta c_n^{-1}\widehat{f}(a_1)} - 1\vert
\exp\{-{{1}\over{2}}
n(1-\widehat{\phi}_1(c_n^{-1}\theta))
\} \nonumber\\
& & \qquad
+(2^{-p} + (6+\widehat{\psi}_1)\widehat{\psi}_{m+1})n{\rm E}\vert
e^{i\theta c_n^{-1}\widehat{f}(a_1)} - 1\vert. \nonumber
\end{eqnarray}
Now, observe that
\begin{equation}
\label{ehei}
\sqrt{n}{\rm E}\vert e^{i\theta c_n^{-1}\widehat{f}(a_1)} - 1\vert
\leq \sqrt{{{n}\over{c_n^2}}}\vert\theta\vert{\rm E}\vert \widehat{f}(a_1) \vert
=o(1).
\end{equation}
By (\ref{ehei})
we can put $m=p\equiv\lfloor\sqrt[4]{n}\rfloor$ in (\ref{eh1}) and (\ref{eh2})
so that by (\ref{ehe})
we get
\[ \overline{\lim}_{\,n}\>
\vert \widehat{\phi}_n(c_n^{-1}\theta) - \exp\{n(\widehat{\phi}_1(c_n^{-1}\theta) -1)\}\vert = 0. \]
Since $\widehat{\phi}_n(c_n^{-1}\theta)\to e^{-\theta^2}$ thus
$n(\widehat{\phi}_1(c_n^{-1}\theta) -1)\to -\theta^2$.
Whence
by the proof of Theorem 8.3.1 in \cite{BGT87} we have that
${\rm E}[\widehat{f}^2(a_1)I_{[\vert \widehat{f}(a_1)\vert\leq x]}]$ is a slowly varying function
in the sense of Karamata and by Lemma \ref{sym} we get that
${\rm E}[f^2(a_1)I_{[\vert f(a_1)\vert\leq x]}]$ varies slowly, too.
Further, by the direct part of this proof we know that
one can choose $d_n=n{\rm E}[f(a_1)]$, which is finite under the slow variation condition.
This completes the proof of Theorem \ref{t2}.

{\bf Acknowledgement.}
I should like to thank
to the referees for
the argument in the line following formula (\ref{esym})
and suggestions that improved
the presentation of the paper.

\end{document}